\documentclass{amsart}
\title{$2$-adic Behavior of Numbers of Domino Tilings}
\newcommand{\eaddress}{\texttt{cohn@math.harvard.edu}}
\author{Henry Cohn\\
Department of Mathematics,
Harvard University\\
Cambridge, MA 02138, USA\\
\eaddress}
\thanks{The author was supported by an NSF Graduate Research
Fellowship.}
\subjclass{05A15, 11A07}
\dedicatory{Dedicated to my grandparents Garnette Cohn (1907--1998)
and Lee Cohn (1908--1998)}
\date{Submitted December 17, 1998; accepted February 18, 1999}

\theoremstyle{plain}
\newtheorem{theorem}{Theorem}

\newtheorem{lemma}[theorem]{Lemma}
\newtheorem{proposition}[theorem]{Proposition}

\newcommand{\Q}{{\mathbb Q}}
\newcommand{\Z}{{\mathbb Z}}

\begin{document}
\thispagestyle{empty}

\begin{abstract}
We study the $2$-adic behavior of the number of domino tilings of a
$2n \times 2n$ square as $n$ varies.  It was previously known that
this number was of the form $2^nf(n)^2$, where $f(n)$ is an odd,
positive integer.  We show that the function $f$ is uniformly
continuous under the $2$-adic
metric, and thus extends to a function on all of $\Z$.  The extension
satisfies the functional equation $f(-1-n) = \pm f(n)$, where the sign
is positive iff $n \equiv 0,3 \pmod{4}$.
\end{abstract}

\maketitle

Kasteleyn \cite{K}, and Temperley and Fisher \cite{TF},
proved that the number of tilings of a $2n \times 2n$
square with $1 \times 2$ dominos is
$$
\prod_{i=1}^n \prod_{j=1}^n \left(4\cos^2\frac{\pi i}{2n+1}
+ 4\cos^2\frac{\pi j}{2n+1}\right).
$$
Although it is by no means obvious at first glance, this number is
always a perfect square or twice a perfect square (see \cite{L}).
Furthermore, it is divisible by $2^n$ but no higher power of $2$.
This fact about $2$-divisibility
was independently proved by several people (see \cite{JSZ},
or see \cite{P} for a combinatorial proof), but there seems to
have been little further investigation of the $2$-adic properties of
these numbers, except for \cite{JS}.

\pagestyle{myheadings}\markboth{\hfill{\sc the electronic journal of
combinatorics 6 (1999), \#R14}}{{\sc the electronic journal of
combinatorics 6 (1999), \#R14}\hfill}

Write the number of tilings as $2^nf(n)^2$, where $f(n)$ is an odd,
positive integer.  In this paper, we study the $2$-adic properties of
the function $f$.  In particular, we will prove the following theorem,
which was conjectured by James Propp:

\begin{theorem}
\label{main}
The function $f$ is uniformly continuous under the $2$-adic metric,
and its unique extension to a function from $\Z_2$ to $\Z_2$
satisfies the functional equation
$$
f(-1-n) = \begin{cases}
f(n) & \textup{if $n \equiv 0,3 \pmod{4}$, and}\\
-f(n) & \textup{if $n \equiv 1,2 \pmod{4}$.}
\end{cases}
$$
\end{theorem}

John and Sachs \cite{JS} have independently investigated the $2$-adic
behavior of $f$, and explicitly determined it modulo $2^6$.
Their methods, as well as ours, can be used to write formulas for $f$
modulo any power of $2$, but no closed form is known.

The proof of Theorem~\ref{main} will not make any use of sophisticated
$2$-adic machinery.
The only non-trivial fact we will require is that the $2$-adic
absolute value extends uniquely to each finite extension of $\Q$.
For this fact, as well as basic definitions and concepts,
the book \cite{G} by Gouv\^ea is an excellent reference.

It is helpful to keep in mind this more elementary description of what
it means for $f$ to be uniformly continuous $2$-adically:  for every
$k$, there exists an $\ell$ such that if $n \equiv m \pmod{2^\ell}$,
then $f(n) \equiv f(m) \pmod{2^k}$.  In particular, we will see that
for our function $f$, the condition $n \equiv m \pmod{2}$ implies that
$f(n) \equiv f(m) \pmod{2}$, and $n \equiv m \pmod{4}$ implies that
$f(n) \equiv f(m) \pmod{4}$.

As a warm-up in using $2$-adic methods, and for the sake of
completeness, we will prove that
that number of tilings of a $2n \times 2n$ square really is of
the form $2^nf(n)^2$, assuming Kasteleyn's theorem.
To do so, we will make use of the fact that the $2$-adic metric
extends to every finite extension of $\Q$, in particular
the cyclotomic extensions, which contain the cosines that appear
in Kasteleyn's product formula.  We can straightforwardly
determine the $2$-adic valuation of each factor, and thus of the
entire product.

Let $\zeta$ be a primitive $(2n+1)$-st root of unity, and define
$$
\alpha_{i,j} = \zeta^i + \zeta^{-i} + \zeta^j + \zeta^{-j}.
$$
Then the number of domino tilings of a $2n \times 2n$ square is
\begin{equation}
\label{kprod}
\prod_{i=1}^n \prod_{j=1}^n (4+\alpha_{i,j}).
\end{equation}
To determine the divisibility by $2$, we look at this number as an
element of $\Q_2(\zeta)$.  Because $2n+1$ is odd, the extension
$\Q_2(\zeta)/\Q_2$ is unramified, so $2$ remains prime in $\Q_2(\zeta)$.
We will use $|\cdot|_2$ to denote the unique extension of the $2$-adic
absolute value to $\Q_2(\zeta)$.

\begin{lemma}
\label{val}
For $1 \le i,j \le n$,
we have
$$
|4+\alpha_{i,j}|_2 =
\begin{cases}
1 & \textup{if $i \ne j$, and}\\
1/2 & \textup{if $i=j$}.
\end{cases}
$$
\end{lemma}

\begin{proof}
The number $4+\alpha_{i,j}$ is an algebraic integer, so its
$2$-adic absolute value is at most $1$.  To determine how much
smaller it is,
first notice that
$$
\alpha_{i,j} = (\zeta^i +
\zeta^j)(\zeta^{i+j}+1)\zeta^{-i}\zeta^{-j}.
$$
In order for $4+\alpha_{i,j}$ to reduce to $0$ modulo $2$, we must
have
$$
\zeta^i \equiv \zeta^{\pm j} \pmod{2}.
$$
However, this is impossible unless $i \equiv \pm j \pmod{2n+1}$,
because $\zeta$ has order $2n+1$ in the residue field.  Since
$1 \le i,j \le n$, the only possibility is $i=j$.

In that case, $4+\alpha_{i,i} = 2(2+\zeta^i+\zeta^{-i})$.
In order to have $|4+\alpha_{i,i}|_2 < 1/2$, the second factor would
need to reduce to $0$.  However, that could happen only if
$\zeta^i \equiv \zeta^{-i} \pmod{2}$, which is impossible.
\end{proof}

By Lemma~\ref{val}, the product \eqref{kprod} is divisible by $2^n$ but not
$2^{n+1}$.  The product of the terms with $i=j$, divided by $2^n$,
is
\begin{equation}
\label{one}
\prod_{i=1}^n (2+\zeta^i +\zeta^{-i}),
\end{equation}
which equals 1, as we can prove by writing
$$
\prod_{i=1}^n (2+\zeta^i +\zeta^{-i})
= \prod_{i=1}^n (1+\zeta^i)(1+\zeta^{-i})
= \prod_{i=1}^n (1+\zeta^i)(1+\zeta^{2n+1-i})
= \prod_{i=1}^{2n}(1+\zeta^i) =1;
$$
the last equality follows from substituting $z=-1$ in
$$
z^{2n+1}-1 = \prod_{i=0}^{2n}(z-\zeta^i).
$$
Thus, the odd factor of the number of tilings of a $2n \times 2n$
square is
$$
f(n)^2 = \prod_{1 \le i < j \le n} (4+\alpha_{i,j})^2.
$$
We are interested in the square root of this quantity, not
the whole odd factor.  The positive square root is
$$
f(n) = \prod_{1 \le i < j \le n} (4+\alpha_{i,j})
$$
(notice that every factor is positive).  It is clearly an integer,
since it is an algebraic integer and is invariant under every
automorphism
of $\Q(\zeta)/\Q$.  Thus, we have shown that the number of tilings
is of the form $2^nf(n)^2$, where $f(n)$ an odd integer.

In determining the $2$-adic behavior of $f$, it seems simplest to
start by examining it modulo $4$.  In that case, we have the formula
$$
f(n) \equiv \prod_{1 \le i < j \le n} \alpha_{i,j} \pmod{4},
$$
and the product appearing in it can actually be evaluated explicitly.

\begin{lemma}
\label{prod}
We have
$$
\prod_{1 \le i < j \le n} \alpha_{i,j} =
\begin{cases}
1 & \textup{if $n \equiv 0,1,3 \pmod{4}$, and}\\
-1 & \textup{if $n \equiv 2 \pmod{4}$.}
\end{cases}
$$
\end{lemma}

\begin{proof}
In this proof, we will write $\zeta^*$ to indicate an unspecified power
of $\zeta$.  Because the product in question is real and the only real
power of $\zeta$ is $1$, we will in several cases be able to see that
factors of $\zeta^*$ equal $1$ without having to count the $\zeta$'s.

Start by observing that
\begin{eqnarray*}
\prod_{1 \le i < j \le n} \alpha_{i,j} &=&
\prod_{i=1}^{n-1} \prod_{j=i+1}^n
(\zeta^{i+j}+1)(\zeta^{i-j}+1)\zeta^{-i} \\
&=& \zeta^{*} \prod_{i=1}^{n-1}\prod_{j=i+1}^n (\zeta^{i+j}+1)
(\zeta^{2n+1+i-j}+1) \\
&=& \zeta^{*} \prod_{i=1}^{n-1}\prod_{s=2i+1}^{2n} (\zeta^s+1).
\end{eqnarray*}
(To prove the last line, check that $i+j$ and $2n+1+i-j$ together run
over the same range as $s$.)

In the factors where $i > n/2$, replace $\zeta^s+1$ with
$\zeta^s(\zeta^{2n+1-s}+1)$.  Now for every $i$, it is easy to check
that
$$
\prod_{s=2i+1}^{2n}(\zeta^s+1)\prod_{s=2(n-i)+1}^{2n}(\zeta^{2n+1-s}+1)
= \prod_{s=1}^{2n}(\zeta^s+1)=1.
$$
When $n$ is odd, pairing $i$ with $n-i$ in this way
takes care of every factor except
for a power of $\zeta$,
which must be real and hence $1$.
Thus, the whole product is $1$ when $n$ is odd, as desired.

In the case when $n$ is even, the pairing between $i$ and $n-i$ leaves
the $i=n/2$ factor unpaired.  The product is thus
\begin{equation}
\label{part}
\zeta^{*} \prod_{s=n+1}^{2n}(\zeta^s+1).
\end{equation}
Notice that
\begin{eqnarray*}
\left(\prod_{s=n+1}^{2n}(1+\zeta^s)\right)^2 &=&
\prod_{s=n+1}^{2n}\zeta^s(1+\zeta^{2n+1-s})
\prod_{s=n+1}^{2n}(1+\zeta^s) \\
&=& \prod_{s=n+1}^{2n} \zeta^s \\
&=& \zeta^{*}.
\end{eqnarray*}
Hence, since every power of $\zeta$ has a square root among the
powers of $\zeta$ (because $2n+1$ is odd),
$$
\prod_{s=n+1}^{2n}(\zeta^s+1) = \pm \zeta^*.
$$
Substituting this result into \eqref{part} shows that the
product we are trying to evaluate must equal $\pm 1$, since
the $\zeta^*$ factor must be real and
therefore $1$.  All that remains is to determine the sign.

Since
$$
\prod_{s=n+1}^{2n} (1+\zeta^s)
$$
and
$$
\prod_{t=1}^n (1+\zeta^t)
$$
are reciprocals, it is enough to answer the question for the second
one (which is notationally slightly simpler).  We know that it is plus
or minus a power of $\zeta$, and need to determine which.  Since $\zeta
= \zeta^{-2n}$, we have
$$
\prod_{t=1}^n(1+\zeta^t) = \prod_{t=1}^n(1+\zeta^{-2nt}) =
\zeta^{*} \prod_{t=1}^n (\zeta^{nt}+\zeta^{-nt}).
$$
The product
$$
\prod_{t=1}^n (\zeta^{nt}+\zeta^{-nt})
$$
is real, so it must be $\pm 1$; to determine which, we
just need to determine its sign.  For that, we write
$$
\zeta^{nt}+\zeta^{-nt} = 2\cos\left(t\pi - \frac{t\pi}{2n+1}\right),
$$
which is negative iff $t$ is odd (assuming $1\le t \le n$).  Thus, the
sign of the
product is
negative iff there are an odd number of odd numbers from $1$ to $n$,
i.e., iff $n \equiv 2 \pmod{4}$ (since $n$ is even).

Therefore, the whole product is $-1$ iff $n \equiv 2 \pmod{4}$, and is
$1$ otherwise.
\end{proof}

Now that we have dealt with the behavior of $f$ modulo $4$, we can
simplify the problem considerably by working with $f^2$ rather than $f$.
Recall that proving uniform continuity is equivalent to showing that
for every $k$, there exists an $\ell$
such that if $n \equiv m \pmod{2^\ell}$, then $f(n) \equiv f(m) \pmod{2^k}$.
If we can find an $\ell$ such that $n \equiv m \pmod{2^\ell}$ implies
that $f(n)^2 \equiv f(m)^2 \pmod{2^{2k}}$, then it follows that
$f(n) \equiv \pm f(m) \pmod{2^k}$, and our knowledge of $f$ modulo $4$
pins down
the sign as $+1$.  The same reasoning applies to the functional equation,
so if we can show that $f^2$ is uniformly continuous $2$-adically and
satisfies $f(-1-n)^2 = f(n)^2$, then we will have proved
Theorem~\ref{main}.

We begin by using \eqref{kprod} to write
\begin{eqnarray*}
2^nf(n)^2 &=& \left(\prod_{i,j = 1}^n \alpha_{i,j}\right)
\prod_{i,j = 1}^n\left(1 + \frac{4}{\alpha_{i,j}} \right) \\
&=& \left(\prod_{i,j = 1}^n \alpha_{i,j}\right) \sum_{k \ge 0}
4^k E_k(n),
\end{eqnarray*}
where $E_k(n)$ is the $k$-th elementary symmetric polynomial in
the $1/\alpha_{i,j}$'s (where $1 \le i,j \le n$).
We can evaluate the product
$$
\prod_{i,j = 1}^n \alpha_{i,j}
$$
by combining Lemma~\ref{prod} with the equation
$$
\prod_{t=1}^n (\zeta^t+\zeta^{-t}) =
(-1)^{\lfloor\frac{n+1}{2}\rfloor},
$$
which can be proved using the techniques of Lemma~\ref{prod}:
it is easily checked that the product squares to $1$, and its sign is
established by writing
$$
\zeta^t+\zeta^{-t}
= 2\cos\frac{2t\pi}{2n+1},
$$
which is positive for $1 \le t < (2n+1)/4$ and negative for $(2n+1)/4
< t \le n$.
This shows that
$$
\prod_{i,j = 1}^n \alpha_{i,j} =
(-1)^{\lfloor\frac{n+1}{2}\rfloor}2^n,
$$
so we conclude that
\begin{equation}
\label{fassum}
f(n)^2 = (-1)^{\lfloor\frac{n+1}{2}\rfloor} \sum_{k \ge 0} 4^k E_k(n).
\end{equation}

The function $n \mapsto (-1)^{\lfloor\frac{n+1}{2}\rfloor}$ is
uniformly continuous 2-adically and invariant under interchanging $n$
with $-1-n$, so to prove these properties for $f^2$ we need only prove
them for the sum on the right of \eqref{fassum}.

Because $\alpha_{i,j}$ has $2$-adic valuation at most $1$,
that of $E_k(n)$ is at least $-k$, and hence $2^kE_k(n)$ is a $2$-adic
integer (in the field $\Q_2(\zeta)$).  Thus, to determine $f(n)^2$
modulo $2^k$ we need only look at the first $k+1$ terms of the sum
\eqref{fassum}.

Define
$$
S_k(n) = \sum_{i,j = 1}^n \frac{1}{\alpha_{i,j}^k}.
$$
We will prove the following proposition about $S_k$.

\begin{proposition}
\label{sk}
For each $k$, $S_k(n)$ is a polynomial over $\Q$ in $n$ and $(-1)^n$.
Furthermore,
$$
S_k(n) = S_k(-1-n).
$$
\end{proposition}

We will call a polynomial in $n$ and $(-1)^n$ a
\textit{quasi-polynomial}.  Notice that every quasi-polynomial over
$\Q$ is uniformly continuous $2$-adically.

In fact, $S_k$ is actually a polynomial of degree $2k$.  However, we
will not need to know that.
The only use we will make of the fact that $S_k$ is a
quasi-polynomial is in proving uniform continuity, so we will prove
only this weaker claim.

Given Proposition~\ref{sk}, the same must hold for $E_k$,
because the $E_k$'s and $S_k$'s are related by the Newton
identities
$$
kE_k = \sum_{i=1}^k (-1)^{i-1} S_iE_{k-i}.
$$
It now follows from \eqref{fassum}
that $f^2$ is indeed uniformly continuous and
satisfies the functional equation.  Thus, we have reduced
Theorem~\ref{main} to Proposition~\ref{sk}.

Define
$$
T_k(n) = \sum_{i,j = 0}^{2n} \frac{1}{\alpha_{i,j}^k},
$$
and
$$
R_k(n) = \sum_{i=0}^{2n} \frac{1}{\alpha_{i,0}^k}.
$$
Because $\alpha_{i,j} = \alpha_{-i,j} = \alpha_{i,-j} = \alpha_{-i,-j}$,
we have
$$
T_k(n) = 4S_k(n) + 2R_k(n) - \frac{1}{\alpha_{0,0}^k}.
$$
To prove Proposition~\ref{sk}, it suffices to prove that
$T_k$ and $R_k$ are quasi-polynomials over $\Q$,
and that $T_k(-1-n) = T_k(n)$ and $R_k(-1-n)
= R_k(n)$.

We can simplify further by reducing $T_k$ to a single sum, as follows.
It is convenient to write everything in terms of roots of unity, so
that
$$
T_k(n) = \sum_{\zeta, \xi} \frac{1}{(\zeta+1/\zeta+\xi+1/\xi)^k},
$$
where $\zeta$ and $\xi$ range over all $(2n+1)$-st roots of unity.
(This notation supersedes our old use of $\zeta$.)
Then we claim that
$$
T_k(n) = \left(\sum_{\zeta} \frac{1}{(\zeta+1/\zeta)^k}\right)^2.
$$
To see this, write the right hand side as
$$
\left(\sum_{\zeta} \frac{1}{(\zeta+1/\zeta)^k}\right)
\left(\sum_{\xi} \frac{1}{(\xi+1/\xi)^k}\right)
=
\sum_{\zeta,\xi} \frac{1}{(\zeta\xi + 1/(\zeta\xi)
+ \zeta/\xi + 1/(\zeta/\xi))^k},
$$
and notice that as $\zeta$ and $\xi$ run over all $(2n+1)$-st roots of
unity, so do $\zeta\xi$ and $\zeta/\xi$.  (This is equivalent
to the fact that every $(2n+1)$-st root of unity has a unique
square root among such roots of unity, because that implies that
the ratio $\xi^2$ between
$\zeta\xi$ and $\zeta/\xi$ does in fact run over all $(2n+1)$-st
roots of unity.)

We can deal with $R_k$ similarly: as $\xi$ runs over all $(2n+1)$-st
roots of unity, so does $\xi^2$, and hence
$$
R_k(n) = \sum_{\zeta} \frac{1}{(2+\zeta+1/\zeta)^k}
= \sum_{\xi} \frac{1}{(2+\xi^2+1/\xi^2)^k}
= \sum_{\xi} \frac{1}{(\xi+1/\xi)^{2k}}.
$$

Define
$$
U_k(n) = \sum_{\zeta} \frac{1}{(\zeta+1/\zeta)^k}.
$$
Now everything comes down to proving the following proposition:

\begin{proposition}
The function $U_k$ is a quasi-polynomial over $\Q$,
and satisfies
$$
U_k(-1-n) = U_k(n).
$$
\end{proposition}

\begin{proof}
The proof is based on the observation that
for any non-zero numbers, the power sums of their
reciprocals are minus the Taylor coefficients of the logarithmic
derivative of the polynomial with those numbers as roots, i.e.,
\begin{eqnarray*}
\frac{d}{dx} \log \prod_{i=1}^m (x-r_i) &=&
\sum_{i=1}^m \frac{1}{x-r_i}\\
&=& \sum_{i=1}^m \frac{-1/r_i}{1-x/r_i}\\
&=& -\sum_{i=1}^m \left(\frac{1}{r_i} + \frac{x}{r_i^2}
+ \frac{x^2}{r_i^3} + \dots\right).
\end{eqnarray*}

To apply this fact to $U_k$, define
\begin{eqnarray*}
P_n(x) &=& \prod_{\zeta}(x-(\zeta+1/\zeta))\\
&=& \prod_{j=0}^{2n} (x-2\cos(2\pi j/(2n+1)))\\
&=& 2(\cos((2n+1)\cos^{-1}(x/2))-1).
\end{eqnarray*}
Then
$$
\frac{d}{dx} \log P_n(x) = \frac{2n+1}{2\sqrt{1-x^2/4}}
\frac{\sin((2n+1)\cos^{-1}(x/2))}{\cos((2n+1)\cos^{-1}(x/2))-1}.
$$
This function is invariant under interchanging $n$ with $-1-n$
(equivalently, interchanging $2n+1$ with $-(2n+1)$),
so its Taylor coefficients are as well.
By the observation above, the coefficient of $x^k$ is $-U_{k+1}(n)$.
Straightforward calculus shows that these coefficients are polynomials
over $\Q$ in $n$, $\sin((2n+1)\pi/2)$, and $\cos((2n+1)\pi/2)$.  Using
the fact that $\cos((2n+1)\pi/2) = 0$ and $\sin((2n+1)\pi/2) = (-1)^n$
completes the proof.
\end{proof}

\section*{Acknowledgements}

I am grateful to James Propp for telling me of his conjecture,
to the anonymous referee for pointing out a sign error in the original
manuscript,
and to Karen Acquista, Noam Elkies, Matthew Emerton, and Vis Taraz for
helpful conversations.


\begin{thebibliography}{JSZ}

\bibitem[G]{G} F.~Gouv\^ea, $p$-adic Numbers: An Introduction, 2nd
ed., Springer-Verlag, New York, 1997.

\bibitem[JS]{JS} P.~John and H.~Sachs, {\sl On a
strange observation in the theory of the dimer problem\/},
preprint, 1998.

\bibitem[JSZ]{JSZ} P.~John, H.~Sachs, and H.~Zernitz, {\sl Problem 5.\
Domino covers in square chessboards\/}, Zastosowania Matematyki
(Applicationes Mathematicae) {\bf XIX} 3--4 (1987), 635--641.

\bibitem[K]{K} P.~W.\ Kasteleyn,
{\sl The statistics of dimers on a lattice, I.
The number of dimer arrangements on a quadratic lattice\/},
Physica {\bf 27} (1961), 1209--1225.

\bibitem[L]{L} L.~Lov\'asz, Combinatorial Problems and Exercises,
North-Holland Publishing Company, Amsterdam, 1979.

\bibitem[P]{P} L.~Pachter, {\sl Combinatorial approaches
and conjectures for $2$-divisibility problems concerning
domino tilings of polyominoes\/}, Electronic Journal of
Combinatorics {\bf 4} (1997), \#R29.

\bibitem[TF]{TF} H.~N.~V.\ Temperley and M.~E.\ Fisher,
{\sl Dimer problem in statistical mechanics---an exact result\/},
Phil.\ Mag. {\bf 6} (1961), 1061--1063.

\end{thebibliography}
\end{document}